\documentclass[12pt,reqno]{amsart}
\usepackage{amssymb,delarray}
\usepackage{amsfonts}
\usepackage{epsfig}
\usepackage[all]{xy}
\usepackage{amscd}
\usepackage{epigraph}

\makeindex{}

\newtheorem{thm}{Theorem}
\newtheorem{lem}{Lemma}
\newtheorem{prop}{Proposition}

\newtheorem{rem}{Remark}

{\catcode`\@=11
\gdef\n@te#1#2{\leavevmode\vadjust{%
 {\setbox\z@\hbox to\z@{\strut#1}%
  \setbox\z@\hbox{\raise\dp\strutbox\box\z@}\ht\z@=\z@\dp\z@=\z@%
  #2\box\z@}}}
\gdef\leftnote#1{\n@te{\hss#1\quad}{}}
\gdef\rightnote#1{\n@te{\quad\kern-\leftskip#1\hss}{\moveright\hsize}}
\gdef\?{\FN@\qumark}
\gdef\qumark{\ifx\next"\DN@"##1"{\leftnote{\rm##1}}\else
 \DN@{\leftnote{\rm??}}\fi{\rm??}\next@}}

\begin{document}

\baselineskip=14.pt plus 2pt 

\title[A Remark on Classical Pluecker's formulae]
{A Remark on Classical Pluecker's formulae}
\author[Vik.S.~Kulikov]{Vik.S. Kulikov}

\address{Steklov Mathematical Institute}  

\email{kulikov@mi.ras.ru}

\thanks{This
research was partially supported by grants of NSh-4713.2010.1, RFBR
08-01-00095, and by AG Laboratory HSE, RF government grant, ag.
11.G34.31.0023. }


\begin{abstract}
For any reduced curve $C\subset \mathbb P^2$, we define the notions
of the number of its virtual cusps $c_v$ and the number of its
virtual nodes $n_v$ which are non-negative,  coincide respectively
with the numbers of ordinary cusps and nodes in the case of cuspidal
curves, and if $\hat C$ is the dual curve of an irreducible curve
$C$ and $\hat n_v$ and $\hat c_v$ are the numbers of its virtual
nodes and virtual cusps, then the integers $c_v$, $n_v$, $\hat c_v$,
$\hat n_v$ satisfy Classical Pl\"{u}cker's formulae.
\end{abstract}

\maketitle

\setcounter{tocdepth}{2}

\def\st{{\sf st}}

\section*{Introduction.}

Let $C\subset \mathbb P^2$ be a reduced curve defined over the field
of complex numbers $\mathbb C$. A curve $C$ is called {\it cuspidal}
if the singular points of $C$ are only the ordinary cusps and nodes.

In modern textbooks on algebraic geometry, classical
Pl$\ddot{u}$cker's formulae are stated as follows (see, for example,
\cite{Br}, \cite{G-H}). \newline {\bf Classical Pl\"{u}cker's
formulae.} {\it Let $C\subset \mathbb P^2$ be an irreducible
cuspidal curve of genus $g$, degree $d\geq 2$, having $c$ ordinary
cusps and $n$ nodes Assume that the dual curve $\hat C$ of $C$ is
also a cuspidal curve. Then
\begin{equation} \label{p1}
\hat d=d(d-1)-3c-2n;
\end{equation}
\begin{equation} \label{p2}
\displaystyle g=\frac{(d-1)(d-2)}{2}-c-n;
\end{equation}
\begin{equation} \label{p3}
d=\hat d(\hat d-1)-3\hat c-2\hat n;
\end{equation}
\begin{equation} \label{p4}
\displaystyle g=\frac{(\hat d-1)(\hat d-2)}{2}-\hat c-\hat n,
\end{equation}
where $\hat c$ and $\hat n$ are the numbers of ordinary cusps and
nodes of $\hat C$ and $\hat d=\deg \hat C$}. \\

Denote by  $V(d,c,n)\subset \mathbb P^{\frac{d(d+3)}{2}}$ the
variety parametrizing the irreducible cuspidal curves of degree $d$
with $c$ ordinary cusps and $n$ nodes.  Very often, if for given
$d$, $c$, and $n$ one of the invariants $\hat c$ or $\hat n$,
obtained as the solution of (\ref{p1}) -- (\ref{p4}), is negative,
then  it is claimed that this is sufficient for the "proof" of the
emptiness of $V(d,c,n)$. However, the correctness of the following
statement is unknown: "{\it the dual curve $\hat C$ of a curve $C$
corresponding to a generic point of $V(d,c,n)$ is cuspidal"}.
Therefore, in general case, it is impossible to conclude the
non-existence of cuspidal curve $C$ if $\hat c$ or $\hat n$ is
negative. Of course, to avoid this problem, it is possible to use
generalized Pl\"{u}cker's formulae including the numbers of all
possible types of singular points of $\hat C$. But, we again have a
difficulty, namely, in this case we must take into account too many
unknown invariables.

To obviate the arising difficulty, in section 1 for any reduced
plane curve $C$ we define the notions of the {\it number} of its
{\it virtual cusps} $c_v$ and the {\it number} of its {\it virtual
nodes} $n_v$ which are non-negative, coincide respectively with the
numbers of ordinary cusps and nodes in the case of cuspidal curves,
and if the dual curve $\hat C$ of an irreducible curve $C$ has $\hat
n_v$ virtual nodes and $\hat c_v$ virtual cusps, then the integers
$c_v$, $n_v$, $\hat c_v$, and $\hat n_v$ satisfy Classical
Pl\"{u}cker's formulae.

In section 2, we investigate the behaviour of the Hessian curve
$H_C$ of a cuspidal curve $C$ at cusps and nodes of $C$, and in
section 3, we give a proof of some inequalities for the numbers of
cusps and nodes of plane cuspidal curves of degree $d$ which was
obtained early in \cite{L} under additional assumption that the dual
curve of a generic cuspidal curve is also cuspidal.

\section{The numbers of virtual cusps and nodes}
Let $(C,p)\subset (\mathbb P^2,p)$ be a germ of a reduced plane
singularity. It splits into several irreducible germs:
$(C,p)=(C_{1},p)\cup\dots\cup (C_{k},p)$. Denote by $m_{j}$ the
multiplicity of the singularity $(C_j,p)$ at the point $p$ and let
$\delta_{p}$ be the $\delta$-invariant of the singularity $(C,p)$.
By definition, the integers
$$c_{v,p}:=\sum_{i=1}^k(m_i-1)$$ and
$$n_{v,p}:=\delta_p-\sum_{i=1}^k(m_i-1)$$ are called respectively the  {\it numbers of virtual cusps}
and {\it virtual nodes} of the singularity $(C,p)$. We have
$\delta_p=c_{v,p}+n_{v,p}$.

\begin{lem} \label{l1} Let $(C,p)\subset (\mathbb P^2,p)$ be a germ of a reduced plane
singularity, $c_{v,p}$ be the number of its virtual cusps and
$n_{v,p}$ be the number of its virtual nodes. Then
\begin{itemize}
\item[($i$)] $c_{v,p}\geq 0$, $n_{v,p}\geq 0$;
\item[($ii$)] if $(C,p)$ is an ordinary cusp, then $c_{v,p}=1$ and
$n_{v,p}=0$;
\item[($iii$)] if $(C,p)$ is an ordinary node, then $c_{v,p}=0$ and
$n_{v,p}=1$.
\end{itemize}
\end{lem}
\proof We prove only the inequality $n_v\geq 0$, since all the other
claims of Lemma \ref{l1} are obvious.  Let
$(C,p)=(C_{1},p)\cup\dots\cup (C_{k},p)$ and $m_i$ be the
multiplicity of its irreducible branch $(C_i,p)$. Then the
multiplicity of $(C,p)$ at $p$ is equal to $m_p=\displaystyle
\sum_{i=1}^k m_i$ and we have
$$\begin{array}{ll} \displaystyle
n_{v,p}= & \delta_p- \displaystyle  \sum_{i=1}^k(m_{i}-1)\geq
\delta_p-
 \sum_{i=1}^km_{i}+1= \\ & \displaystyle \delta_p- (m_p-1)\geq
 \delta_p-\frac{m_p(m_p-1)}{2}\geq 0,\end{array}$$
since $m_p\geq 2$ for singular points and $\displaystyle
\delta_p\geq \frac{m_p(m_p-1)}{2}$. Therefore, we have $$
\displaystyle n_v=\sum_{p\in \text{Sing}\, C}n_{v,p}\geq 0. \qed $$
\\

Let $C\subset \mathbb P^2$ be a reduced curve. Denote by
$\text{Sing}\, C$ the set of its singular points. By definition, we
put
$$c_v:=\sum_{p\in \text{Sing}\, C}c_{v,p},$$
$$n_v:=\sum_{p\in \text{Sing}\, C}n_{v,p}$$
and call these integers respectively the  {\it number of virtual
cusps} and the {\it number virtual nodes} of the curve $C$. If $C$
is an irreducible curve of degree $d$ and geometric genus $g$, then
we have $g=\displaystyle \frac{(d-1)(d-2)}{2}-\delta_C$, where
$\delta_C=\displaystyle \sum_{p\in\text{Sing}\, C}\delta_p$ is the
$\delta$-invariant of $C$. Therefore, we have

\begin{equation} \label{genus} g=\frac{(d-1)(d-2)}{2}-c_v-n_v.\end{equation}

The following proposition is a corollary of Lemma \ref{l1}.
\begin{prop} \label{cl1} Let $c_v$ be the number of virtual cusps and $n_v$ be the number
of virtual nodes of a reduced curve $C\subset \mathbb P^2$. We have
\begin{itemize}
\item[($i$)] $c_v\geq 0$ and $n_v\geq 0$,
\item[($ii$)]
 if $C$ is a cuspidal curve, then
$c_v$ and $n_v$ are equal respectively to the number $c$ of cusps
and the number $n$ of nodes of $C$.
\end{itemize}
\end{prop}
\begin{thm} {\rm (}{\bf Pl\"{u}cker's formulae}{\rm )}. Let $C$ and $\hat{C}$ be
irreducible dual curves of genus $g$, $\deg C=d\geq 2$, $\deg \hat
C=\hat d$, and $c_v$, $n_v$, $\hat c_v$, $\hat n_v$ are the numbers
of their virtual cusps and nodes, respectively. Then we have the
following equalities:
\begin{equation} \label{p11}
\hat d=d(d-1)-3c_v-2n_v;
\end{equation}
\begin{equation} \label{p21}
\displaystyle 2g=(d-1)(d-2)-2c_v-2n_v;
\end{equation}
\begin{equation} \label{p31}
d=\hat d(\hat d-1)-3\hat c_v-2\hat n_v;
\end{equation}
\begin{equation} \label{p41}
\displaystyle 2g=(\hat d-1)(\hat d-2)-2\hat c_v-2\hat n_v.
\end{equation}
\end{thm}
\proof To prove Pl\"{u}cker's formulae, we need the following
\begin{lem} For an irreducible plane curve
$C$ we have
\begin{equation}\label{d^*} \hat d= 2d+2(g-1)-c_v,
\end{equation}
\begin{equation}\label{c^*} \hat c_v= 3d+6(g-1)-2c_v,
\end{equation}
\begin{equation}\label{d} d= 2\hat d+2(g-1)-\hat c_v,
\end{equation}
\begin{equation}\label{c} c_v= 3\hat d+6(g-1)-2\hat c_v.
\end{equation}
\end{lem}
\proof Denote by $\nu: \overline C\to C$ and $\hat{\nu}: \overline
C\to \hat C$ the normalization morphisms, consider generic (with
respect to $C$ and $\hat C$) linear projections $\text{pr}:\mathbb
P^2\to\mathbb P^1$ and $\hat{\text{pr}}:\hat{\mathbb P}^2\to\mathbb
P^1$, and put $\pi=\text{pr}\circ \nu$ and
$\hat{\pi}=\hat{\text{pr}}\circ \hat{\nu}$. We have $\deg\pi=d$ and
$\deg\hat{\pi}=\hat{d}$.

Let $\nu^{-1}(x_i)=\{ y_{i,1},\dots ,y_{i,m_i}\}$ for $x_i\in
\text{Sing}\, C$. For each point $y_{i,j}$ denote by $r_{i,j}$ the
ramification index of $\pi$ at $y_{i,j}$. It is easy to see that
$r_{i,j}$ coincides with the multiplicity $m_{i,j}$ at $x_i$ of the
irreducible germ $(C_{i,j},x_i)\subset (C,x_i)$ corresponding to the
point $y_{i,j}$. Therefore, we have
$$c_v=\sum_{i,j}(r_{i,j}-1).$$

Applying Hurwitz formula to $\pi$ and $\hat{\pi}$, we obtain
\begin{equation}
2(g-1)=-2d+c_v+\hat d\end{equation} and
\begin{equation} \label{p_a}
2(g-1)=-2\hat d+\hat c_v+d\end{equation} which give formulae
(\ref{d^*}) and (\ref{d}).

To prove (\ref{c^*}), note that $\hat c_v=2\hat d+2(g-1)-d$ by
(\ref{d}). Therefore $$\hat c_v=2(2d+2(g-1)-c_v) +2(g-1)-d$$ by
(\ref{d^*}), that is, $\hat c_v= 3d+6(g-1)-2c_v$. Formula (\ref{c})
is obtained similarly. \qed \\

It follows from (\ref{genus}) that
\begin{equation} \label{argen}
2(g-1)+2c_v+2n_v=d(d-3), \qquad 2(g-1)+2\hat c_v+2\hat n_v=\hat
d(\hat d-3)\end{equation} which are equivalent to (\ref{p21}) and
(\ref{p41}). To complete the proof of Pl\"{u}cker's formulae, notice
that formulae (\ref{p11}) and (\ref{p31}) easily follow from
equations (\ref{d^*}) -- (\ref{c}) and (\ref{argen}). \qed

\section{On the Hessian curve of a cuspidal curve}
Let $C\subset \mathbb P^2$ be an irreducible cuspidal curve of
degree $d$ with $c$ cusps and $n$ nodes. It follows from (\ref{p21})
and (\ref{c^*}) that
\begin{equation}\label{eqmain} 8c+6n +\hat c_v= 3d(d-2).\end{equation}

Equality (\ref{eqmain}) has a natural geometric meaning. To explain
it, let the curve $C$ is given by equation $F(x_0,x_1,x_2)=0$, where
$x_0,x_1,x_2$ are homogeneous coordinates in $\mathbb P^2$. Consider
the Hessian curve $H_C\subset \mathbb P^2$ of the curve $C$. It is
given by equation $\det (\frac{\partial^2F}{\partial x_i\partial
x_j})=0$. We have $\deg H_C=3(d-2)$. Therefore the intersection
number $(C,H_C)_{\mathbb P^2}$ is equal to $3d(d-2)$. On the other
hand, it is well-known (see, for example, \cite{Br}) that the curves
$C$ and $H_C$ meet at the singular points and at the inflection
points of the curve $C$. Therefore we have
\begin{equation} \label{int} \Sigma^{\prime}(C,H_C)_{p}+\Sigma^{\prime\prime}(C,H_C)_{p}+
\Sigma^{\prime\prime\prime}(C,H_C)_{p}=(C,H_C)_{\mathbb
P^2}=3d(d-2),\end{equation} where $(C,H_C)_p$ is the intersection
number of the curves $C$ and $H_C$ at a point $p\in C$ and the sum
$\sum^{\prime}$ is taken over all cusps of $C$, the sum
$\sum^{\prime\prime}$ is taken over all nodes of $C$, and the sum
$\sum^{\prime\prime\prime}$ is taken over all inflection points of
$C$.

Let us show that the coefficients involving in equation
(\ref{eqmain}) have the following geometric meaning: {\it equality}
(\ref{eqmain}) {\it is the same as equality} (\ref{int}), that is,
the coefficient $8$ in (\ref{eqmain}) is the intersection number
$(C,H_C)_p$ at a cusp $p\in C$, the coefficient $6$ is the
intersection number $(C,H_C)_p$ at a node $p\in C$, and $\hat
c_v=\sum^{\prime\prime\prime}(C,H_C)_{p}$. Indeed, let $p$ be a cusp
of $C$. Without loss of generality, we can assume that $p=(0,0,1)$
and
$$F(x_0,x_1,x_2)=x_0^2U(x_0,x_1,x_2)+
x_0x_1^2V(x_0,x_1,x_2)+x_1^3W(x_0,x_1,x_2),$$ where $U$ is a
homogeneous polynomial of degree $d-2$ such that $U(0,0,1)=1$ and
$V$ and $W$ are homogeneous polynomials of degree $d-3$ such that
$W(0,0,1)=1$. Put $a=V(0,0,1)$, then in non-homogeneous coordinates
$x=\frac{x_0}{x_2}, y=\frac{x_1}{x_2}$ we have $p=(0,0)$, the curve
$C$ is given by equation of the form
$$x^2+y^3+axy^2+bx^2y+cx^3+ \text {terms of higher degree}=0,$$ and
the curve $H_C$  is given by equation of the form
$$x^2(6y+2ax)+ \text {terms of higher degree}=0.$$ Easy computation
(applying $\sigma$-process with center at $p$) gives the following
inequality:
\begin{equation}\label{1} (C,H_C)_p\geq 8\end{equation} if $p$ is a
cusp of $C$.

Let $p$ be a node of $C$. Again, without loss of generality, we can
assume that $p=(0,0,1)$ and
$$F(x_0,x_1,x_2)=x_0x_1U(x_0,x_1,x_2)+
V(x_0,x_1)W(x_0,x_1,x_2),$$ where $U$ is a homogeneous polynomial of
degree $d-2$ such that $U(0,0,1)=1$, $V$ is a homogeneous polynomial
of degree $3$, and $W$ is a homogeneous polynomial of degree $d-3$.
In non-homogeneous coordinates $x=\frac{x_0}{x_2},
y=\frac{x_1}{x_2}$ we have $p=(0,0)$, the curve $C$ is given by
equation of the form
$$xy+ \text {terms of higher degree}=0,$$ and
the curve $H_C$  is given by equation of the same form
$$xy+ \text {terms of higher degree}=0.$$ Easy computation
(applying $\sigma$-process with center at $p$) gives the following
inequality:   \begin{equation}\label{2} (C,H_C)_p\geq
6\end{equation} if $p$ is a node of $C$.

If $p$ is an $r$-tuple inflection point of $C$ (that is,
$(C,L_p)_p=r+2$, where the line $L_p$ is tangent to $C$ at $p$),
then by Theorem 1 on page 289 in \cite{Br}, we have $(C,H_C)_p=r$.
On the other hand, the branch $(\hat C_i,\hat p)$ of the dual curve
$\hat C$, corresponding to an irreducible branch $(C_i,p)\subset
(C,p)$ at a point $p$ of a cuspidal curve $C$, is singular if and
only if $p$ is an inflection point of $C$; and the branch $(\hat
C,\hat p)$, corresponding to the branch $(C,p)$ at $r$-tuple
inflection point $p\in C$, has a singularity of type
$u^{r+1}-v^{r+2}=0$. The multiplicity $m_{\hat p}$ of this
singularity is equal to $r+1$. Therefore, we have
\begin{equation} \label{3} \Sigma^{\prime\prime\prime}(C,H_C)_{p}=
\sum_{(\hat C,\hat p)} (m_{\hat p}-1) =\hat c_v.\end{equation}

Finally, it follows from (\ref{eqmain}) -- (\ref{3}) that
inequalities (\ref{1}) and (\ref{2}) are the equalities in the case
of cuspidal curves.

\section{Lefschetz's inequalities}
As above, let $C\subset \mathbb P^2$ be an irreducible cuspidal
curve of degree $d$ and genus $g$ having $c$ cusps and $n$ nodes.

In \cite{L},  assuming that for a generic cuspidal curve with given
numerical invariants the dual curve is also cuspidal, Lefschetz
proved the following inequalities
\begin{equation} \label{ineven} c\leq \frac{3}{2}d+3(g-1)
\end{equation} if $d$ is even and
\begin{equation} \label{inodd} c\leq \frac{3d-1}{2}+3(g-1)
\end{equation} if $d$ is odd. It follows from  (\ref{c^*}) that these inequalities
occur for any plane cuspidal curve, since $\hat c_v$ is a
non-negative integer.

Note also that equality (\ref{eqmain}) gives rise to the following
statement: {\it for a plane cuspidal curve of degree $d\geq 2$ the
following inequality holds}: \begin{equation} \label{inmain}
8c+6n\leq 3d(d-2)-\frac{1-(-1)^d}{2}.\end{equation}

\begin{rem}{\rm One can show that} for any $d=2k$, $k\geq 3$, and
for any $g\geq 0$ such that $2\leq 3g\leq k-4$ or $g\leq 1$, there
exist a cuspidal curve of degree $d$ having $c=3(k+g-1)$ cusps and
$n=2(k-1)(k-2)-4g$ nodes for which inequality {\rm (\ref{inmain})}
becomes the equality. If $d=2k+1$, $k\geq 3$, then for any $g$ such
that $2\leq 3g\leq k-4$ or $g\leq 1$, there exist a cuspidal curve
of degree $d$ having $c=3(k+g)-2$ cusps, $n=2(k-1)^2-4g$ nodes, and
for which inequality {\rm (\ref{inmain})} becomes the equality. {\rm
The proof of these statements follows from the fact that the genus
of such curves $C$ is equal to $g$ and for these curves the dual
curves $\hat C$ have degree $\hat d=2(g-1)+7 + \frac{1-(-1)^d}{2}$
and the number of virtual cusps $\hat c_v=\frac{1-(-1)^d}{2}$.
Therefore such curves can be obtained as the image of a generic
linear projection to $\mathbb P^2$ of a smooth curve $\overline
C\subset \mathbb P^{\hat d-g}$ of degree $\hat d$ birationally
isomorphic to $C$. Standard computations (which we leave to the
reader) of codimension of the locus of "bad" projections shows that
in this case there is a linear projection $\text{pr}: \mathbb
P^{\hat d-g}\to \mathbb P^2$ such that $\text{pr}(\overline C)=\hat
C$ is a cuspidal curve with $\hat c=\frac{1-(-1)^d}{2}$ and its dual
curve $C$ is also cuspidal}.
\end{rem}

For completeness, let me remind also the following well-known
inequalities which we have for plane cuspidal curves:
\begin{equation} 3c+2n< d(d-1)-\sqrt{d},\end{equation} \begin{equation} 2c+2n\leq (d-1)(d-2),\end{equation}
\begin{equation} \label{fina} d(d-2)(d^2-9)+(3c+2n)^2+27c+20n\geq 2d(d-1)(3c+2n)\end{equation} which are
consequences of equalities (\ref{p11}) -- (\ref{p41}) and the
inequalities  $\hat d\geq \sqrt d$, $g\geq 0$,  $\hat n_v\geq 0$;
$$16c+9n\leq d(5d-6)$$ ({\it Hirzebruch -- Ivinskis
inequality} (\cite{Hirz})) which is true in the case of even $d$.

Note also that inequalities (\ref{inmain}) -- (\ref{fina}) hold for
any irreducible plane curve if we substitute $c_v$ and $n_v$ instead
of $c$ and $n$.

\end{document}